\documentclass[12pt]{amsart}
\usepackage{amsmath}
\usepackage{amssymb}
\newcommand{\newpars}{\par\smallskip\noindent}

\newcommand{\newparb}{\par\bigskip\noindent}

\newcommand{\tmod}[1]{\allowbreak\ ({\textrm{mod}}\,\,#1)}
\newcommand{\hMc}{\Mc}
\newcommand{\sMc}{M\protect\raisebox{3pt}{\mbox{\protect\scriptsize c}}}
\newcommand{\tMc}{M\protect\raisebox{3pt}{\mbox{\protect\tiny c}}}

\newcommand{\Fano}{\operatorname{Fano}}
\newcommand{\Mathieu}{\operatorname{Mathieu}}
\newcommand{\Serre}{\operatorname{Serre}}

\newtheorem{thm}{Theorem}[section]

\newtheorem{prop}[thm]{Proposition}
\theoremstyle{definition}

\theoremstyle{remark}

\begin{document}
\title{Galois groups over function fields of positive characteristic}
\author[Conway]{John Conway}
\author[McKay]{John {\hMc}Kay}
\author[Trojan]{Allan Trojan}
\address{John Conway,
         Department of Mathematics, Fine Hall,
         Princeton University, Washington Road,
         Princeton, N.J. 08544-1000  USA}
         \email{conway@math.princeton.edu}
\address{John {\tMc}Kay,
         Department of Mathematics and CICMA
         Concordia University
         1455 de Maisonneuve Blvd. West, Montreal, Quebec H3G 1M8, Canada.}
         \email{mac@mathstat.concordia.ca}
\address{Allan Trojan, SASIT, Atkinson College,
         York University,
         4700 Keele St.,
         Toronto, Ontario M3J 1P3, Canada}
         \email{atrojan@yorku.ca}

\begin{abstract}
We describe examples motivated by the work of Serre and Abhyankar.
\end{abstract}
\maketitle
\par
\section{The main result}
Let $K$ be a field of characteristic $p$ and $\overline{K}$ its
algebraic closure, $K(t)$ the field of functions in variable $t$,
and $q$ a power of $p$; Galois fields of
order $q$ will be denoted by $F_q$.
A survey of computational Galois theory is found in the special issue 
of J. Symbolic Computation \cite{JSC}; here we describe techniques for 
computing Galois groups / $K(t)$ in maple.

Certain questions and conjectures concerning Galois 
theory over $\overline{K}(t)$ were raised by Abhyankar
\cite{Abh1} in 1957.  Apparently the first result was 
obtained in 1988 by Serre (in Abhyankar,\cite{Abh2},appendix),
who proved that $PSL_2(q)$ occurs for the
polynomial $x^{q+1}-tx+1$.

Abhyankar continued in \cite{Abh1}, obtaining results
for unramified coverings of the form:
\[
x^{n}-at^{u}x^{v}+1,
\quad
(v,p)=1,
\quad
n=p+v,
\]
and
\[
x^{n}-ax^{v}+t^{u},
\quad
(v,p)=1,
\quad
n\equiv0\tmod p,
\quad
u\equiv0\tmod v.
\]
The groups obtained have the form $S_{n}$,
$A_{n}$, $PSL_2(p)$ or $PSL_2(2^3)$.

They used algebraic geometry to construct a
Galois covering.
Abhyankar used a method that relied on a
characterization
of the Galois groups as permutation groups while
Serre used a method based on L\"uroth's
theorem and the
invariant theory of Dickson.
Later, in \cite{Abh3}, Abhyankar
obtained the
Mathieu group, $M_{23}$, as the Galois group of
$x^{23}+tx^{3}+1\;/\;F_2(t)$.
His proof was based on an idea used by Serre (in
Abhyankar and Yie,
\cite{AbhYie} ) who showed that 
$PSL(3,2)$ is the
Galois group of $x^{7}+tx^{3}+1\;/\;F_2(t).$

They express the irreducible polynomial
as a factor
of an additive polynomial, 
that is, a polynomial of the form:\
$g(x)=\sum_{i=0}^{n}a_{i}x^{p^{i}}$.\

Such a polynomial has the property that its zeros form
a vector space
over $F_p$. Thus the Galois group of any polynomial
dividing $g$
is bounded above by $GL(n,p)$. The action of the Frobenius 
element, and further knowledge on transitivity, 
may provide lower bounds, sometimes determining the group.
 Abhyankar obtained these transitivity
conditions by the technique of 'XXXXthrowing away a root'.
Abhyankar continued these techniques in a long series
of papers. In \cite{Abh4} he has written a expository
summary with many references.

Using simpler arguments, avoiding algebraic geometry, 
and, furthermore, not requiring algebraic closure, we 
\cite{ConMcK} found the Galois groups of several 
polynomials over $F_p(t)$ for $p=2,3$.

If $t$ is algebraic over $F_p$, and the polynomial 
$f(x,t)$ over $F_p(t)$ has distinct roots,
and has irreducible factors of degree 
$d_{1},\,d_{2},\,\ldots$
then the Galois group of $f(x,t)$ contains a permutation of shape 
$d_{1},\,d_{2},\,\ldots$
(van der Waerden, \cite{vderW}).
\par\smallskip
The following Galois groups were obtained over $F_2(t)$:
\[
\begin{tabular}{|c|c|}
\hline
               &            \\[-3.0\jot]
$x^{24}+x+t$   & $M_{24}$   \\
$x^{23}+x^3+t$ & $M_{23}$   \\
$x^{13}+x+t$   & $PSL_3(3)$ \\
$x^8+x^7+t$    & $PSL_2(7)$ \\
$x^7+x+t$      & $PSL_3(2)$ \\
$x^4+x+t$      & $D_{4}$    \\[+0.5\jot]
\hline
\end{tabular}
\]
and over $F_3(t)$:
\[
\begin{tabular}{|c|c|}
\hline
                  &          \\[-3.0\jot]
$x^{12}+x+t$      & $M_{11}$ \\
$x^{11}+tx^{2}-1$ & $M_{11}$ \\[+0.5\jot]
\hline
\end{tabular}
\]
\newparb
Here is the proof for the polynomial $f = x^{24}+x+t$:
\par\smallskip
1. Specialize $t$ to be a root of one of the following equations:
\[ \text{$t=0,t=1,\,t^{3}+t^{2}+1=0$.} \]
2. Over the algebraic closure of $F_2$ we find irreducible factors 
   of respective degree shapes:
\[ \text{$[11^{2}1^{2}]$,$[21,3]$ and $[23,1]$.} \]
3. It follows that the Galois group $G$ (over $F_2(t)$) 
   is at least $3$-transitive and has order a multiple of 7.11.23,
   thus, as a permutation group, it is one of 
   $S_{24}$, $A_{24}$, or $M_{24}$. 

4. Obtain an additive polynomial as a multiple of $f$ as follows:
\begin{align*}
0
&\equiv x^{32}\!+x^{9}\!+t x^{8} \tmod f
\\
&\equiv x^{256}\!+x^{72}\!+t^{8} x^{64} 
&&\text{(raising to the $8^{th}$ power)}
\\
&\equiv x^{256}\!+t^{8} x^{64}\!+x^{3}\!+t x^{2}\!+t^{2} x\!+t^{3}
&&\text{(substituting $x^{72}\equiv(x+t)^{3}$)}
\\
&\equiv 
 x^{2048}\!+t^{64} x^{512}\!+x^{24}\!+t^{8} x^{16}\!+t^{16} x^{8}\!+t^{24}
&&\text{(raising to the $8^{th}$ power)}
\end{align*}
Thus, since $x^{24}=x+t$, we get the semi-linear
relationship:
\begin{equation} \label{semilinear}
t+t^{24} = L = x^{2048}+t^{64}\ x^{512}+t^{8}\
x^{16}+t^{16}\
x^{8}+x
\end{equation}

So \textbf{$L^{2}=(t+t^{24})\ L$} is an equation of
degree $2^{12}$ in $x$ whose terms involve only the
powers of\[
x^{2^{n}}:0\leq n\leq12\]

5. The zeros of $f$ satisfy an additive polynomial and
so lie in a vector
space of dimension 12 over $F_2$. But neither
$A_{24}$ nor $S_{24}$
can act on such a space. $(G \subseteq GL_{12}(2))$ so
$G=M_{24}.$
\newparb
Except for the $M_{11}$ degree 11 case, similar methods are used to 
determine the other Galois groups,
and a computer is used only for the factorization of the polynomial
over finite fields. 
In the last case: $f(x) = x^{11}+ tx^{2}-1$ over $F_3$.
We seek an additive polynomial of degree $3^5$ but find only
an additive polynomial of degree $3^{10}$ with 326 terms. It is
unclear to us whether one can derive a degree $3^5$ additive 
polynomial using this method but modifying the initial $f(x)$.
\newparb
The zeros of $f$ satisfy an additive polynomial
whose zeros form a vector space of dimension 10 over $F_3$.
As before, by factoring $f$ over various finite extensions 
of $F_3(t)$,
G = Gal(f) is shown to contain permutations of shapes:
\
$[5^{2}1],[6,3,2],[8,2,1],[11]$
\
which are shapes occurring in  $M_{11}$.
So the group is at least 3-transitive
which means it is one of the following:\
$S_{11}$ , $A_{11}$ , $M_{11}$
but $S_{11}$ is excluded
since the discriminant of $f(x,t)$ with respect to $x$ is a square, 1.
\par\smallskip
To exclude $A_{11}$ requires more.
The 5-set resolvent polynomial $f^{\times{5}}$ \cite{CasMcK}
is computed over its coefficient ring with a call to 
``galois/rsetpol"$(f,x,5)$.
(The zeros of this polynomial, of degree $\binom{11}{5}=462$,
are the products of the 5-element subsets of the zeros of $f$).
This resolvent polynomial has a factor of $x$-degree 66:
\begin{equation} \label{f66}
\begin{split}
f_{66}=\, &
{x}^{66}+t{x}^{62}+2\,{t}^{5}{x}^{57}+{t}^{3}{x}^{54}+{t}^{6}{x}^{53}+t{x}^{51}+2\,{t}^{7}{x}^{49}+{t}^{10}{x}^{48}+
2\,{t}^{2}{x}^{47}+\\
& {t}^{5}{x}^{46}+  2\,{t}^{8}{x}^{45}+{t}^{3}{x}^{43}+2\,{t}^{6}{x}^{42}+ t{x}^{40}+
2\,{t}^{7}{x}^{38}+ {t}^{2}{x}^{36}+ {t}^{8}{x}^{34}+ {x}^{33}+\\
& 2\,{t}^{6}{x}^{31}+ {t}^{9}{x}^{30}+2\,t{x}^{29}+2\,{t}^{7}{x}^{27}+
2\,{t}^{5}{x}^{24}+2\,{t}^{3}{x}^{21}+2\,{t}^{6}{x}^{20}+2\,t{x}^{18}+\\
& 2\,{t}^{7}{x}^{16}+
{t}^{2}{x}^{14}+{t}^{5}{x}^{13}+
{t}^{8}{x}^{12}+
 2\,{t}^{3}{x}^{10}+2\,{t}^{6}{x}^{9}+2\,t{x}^{7}+2\,{t}^{2}{x}^{3}+1
\end{split}
\end{equation}

so finally the Galois group is $M_{11}$, not $A_{11}$.

\section{Series Solutions to the Equations}

The methods so far are non-constructive. To construct the group
and to study the Galois correspondence, it is useful to
work with exact zeros of the polynomial,
either as formal Taylor series or formal Puiseux series in the variable
$t$ over the
field $F_q$.

We write $\Fano(x)$ for the polynomial $x^{7}+tx+1$ / $F_2(t)$
which we shall see is related closely to Serre's $x^{7}+tx^{3}+1$.

\begin{prop} \label{PFano}
$\Fano(x)$ has 7 formal Taylor series solutions, of the
form:\[
f_{\varepsilon}=\sum_{i\in\mathbf{\Omega}}\varepsilon^{2^{N(i)}}t^{i}\]

where $\mathbf{\Omega}$ is the set of non-negative
integers, $n$,
such that $n=0$ or the lengths of the blocks of zeros
in the binary
expansion of $n$ are all multiples of 3; $N(i)$ is the
number of 1's
in the binary expansion of $i$; and $\varepsilon$ is a
primitive
$7^{th}$ root of unity in $F_{2^3}$.
\end{prop}

We show the series $f_{\varepsilon}$ are solutions to the additive
polynomial $x^{8}+tx^{2}+x$ and we
identify them
with the set of pairs $\Lambda = \{ i\in \Omega :(i,N(i) \textrm{ mod } 3)  \}$.

The map $\Psi_{0}$ that sends the series $f$ to $f^{8}$ sends $(i,N(i))$ to $(j,N(i))$
where the binary expansion of $j$ is obtained
from that of $i$ 
by appending the block 000;
however, $N(i) = N(j)$, so $(i,N(i))$ is sent to $(j,N(j))$ and
so $\Psi_{0}$ sends $\Lambda$ onto the subset $\Lambda_{0}= \{i\in \Omega \textrm{ and } i \equiv
0 \textrm{ mod } 8:(i,N(i))\}$

The map $\Psi_{1}$  that sends the series $f$ to $t f^{2}$ sends $(i,N(i))$ to $(j,N(i)+1)$
where $j$ is obtained from $i$ by appending the digit 1
to the binary representation of $i$ but in this case also, 
$(i,N(i))$ is again sent to $(j,N(j))$
since $N(i)+1 =N(j)$,
so $\Psi_{1}$ sends $\Lambda$ onto the subset 
$\Lambda_{1}= \{i\in \Omega 
              \textrm{ and } i \equiv 1 \textrm{ mod 2 }:(i,N(i))\}$
hence the map $f^{8} + t f^{2} = \Psi_{0}+ \Psi_{1}$ 
is a bijection of $\Lambda$ with $\Lambda$.
\newpars
The preceding proof indicates how to obtain series solutions to an
additive polynomial, however, not every additive polynomial has a 
complete set of solutions in the form of a Taylor series or even a 
Puiseux series.

The polynomial $x^{4}+x^{2}+t x$, for example, has  only one series
solution $(t+t^{3}+t^{5}+t^{9}+...)$ but if we substitute
$t=\frac{1}{s}$ the polynomial equation becomes:
$s x^{4}+s x^{2}+ x$,
which has 3 solutions:

$x=A s^{-\frac{1}{3}}+B\cdot s^{\frac{1}{3}}+A\cdot
s^{\frac{5}{3}}+B\cdot s^{\frac{7}{3}}+...$

where $A$ is any of the three cube roots of unity and
$B=A^{2}$.
These solutions may be thought of as Puiseux series
``about the point at infinity".

When the polynomial is not additive, one can get series solutions by
finding an additive polynomial for which it is a factor.

\begin{prop} \label{PMathieu}
The polynomial $\Mathieu(x)=x^{24}+x+t$ has 23 formal Taylor
series solutions of
the form:\[
f_{\alpha}=\sum_{i\in\Omega}\alpha^{1-i}t^{i}\]

Where $\alpha$ is any non-zero solution to the
equation $x^{24}=x$
in an extension of $F_2$ and $\Omega$ is the set of
non-negative
integers, $n$, such that the binary representation of
$n$ has the
form

0 or 11 or 101000 or 100000000

followed by any combination of

010, 0001, 000001000, 00000000000

followed by 000.
\end{prop}

First, suppose\[
f=\sum_{i\in\Lambda}\alpha^{K(i)}t^{i}\]

is a Taylor series solution to $\Mathieu(x)=0$,

then for every $i$, we have
$K(i)+i\equiv1\pmod{23}$.

This follows by induction from the fact that the terms
of $f^{24}$
are of the form:
$\alpha^{24K(i)}t^{24i}$\,
or\,$\alpha^{16K(i)+8K(i')}t^{16i+8i'}$ so 
such a Taylor
series solution is 
determined by
$\alpha$, which satisfies $\alpha^{23}=1$, and the
exponents of
$t$, which are independent of $\alpha$. To determine
$f(x)$ it is
sufficient to determine the exponents of $t$
occurring, so we may
assume $\alpha=1$.

We make the substitutions, $g=f+t,u=t^{8}$,
in the semi-linear equation (\ref{semilinear})
satisfied by $f$, to get:
\begin{equation*}
g^{2048}+u^{8}g^{512}+u g^{16}+u^{2}g^{8}+g+u^{256}+u^{72}+u^{3}=0.
\end{equation*}
Define the following binary sequences:

$Z=0$, $A=100000000$, $B=1001000$, $C=11$,
$K=000000000000$,\\
$L=000001000$, $M=0001$, $N=010$

For $g$ a Taylor series in $u$ over $F_2$ with constant term
equal to 1, let $\Omega_{g}$ be the set of binary strings which
occur as exponents of $u$ in $g$. Since the constant
term of $g$ is
1 and the terms $u^{256}+u^{72}+u^{3}$ occur in the
semi-linear
polynomial, $Z$, $A$, $B$, $C\in$ $\Omega_{g}$.

As in the proof of Proposition~\ref{PFano}  we can
consider
$g^{2048}$, $u^{8} \cdot g^{512}$, etc. as operators
on the set of
binary strings. $g^{2048}$, $u^{8} \cdot g^{512}$
\ldots append the
string $K$, $L$ \ldots respectively. Since the images
of these four
operators are disjoint, if the string  $S\in
{\Omega_{g}}$,then $S
\circ K$, $S \circ L$, $S \circ M$, $S \circ N \in
\Omega_{g}$,
where $ \circ $ represents string concatenation.
\newpars
Remark: The remaining series solution, with constant term 0, can be
found using the same conditions omitting the initial term $Z$.
\newpars
We return now to the polynomial $\Fano(x)$ of
Proposition~\ref{PFano}
and its 'complementary' polynomial $\Serre(x)$.

It is instructive to see how the Galois group can be determined just
from the series solutions determined by Proposition~\ref{PFano}
as well as the fields given by the Galois correspondence.

In the notation of Proposition~\ref{PFano}, let $\varepsilon$ satisfy 
$\varepsilon^{2}+\varepsilon+1=0$ and
\begin{align*}
&P_{i}=f_{\varepsilon^{i}},
\\
&L_{0}=P_{1} \cdot P_{2} \cdot P_{4},\,
 L_{1}=P_{2} \cdot P_{3} \cdot P_{5},\,
 \ldots,\,
 L_{6}=P_{0} \cdot P_{1} \cdot P_{3}
\\
&\text{(corresponding to the lines in the Fano projective plane)},
\\
&K=F_2(t,P_{0},P_{1},...,P_{6}),
\end{align*}
and let $R$ be the subfield of $K$ containing all
Taylor series in
$t$ whose coefficients are in $F_2$.
\begin{prop}.

\begin{enumerate}
\item
\begin{tabbing}
$P_{0},P_{1},...,P_{6}$ are the zeros of
$\Fano(x)=x^{7}+t\cdot x+1$
\\ $L_{0},L_{1},...,L_{6}$ are the zeros of
$\Serre(x)=x^{7}+s\cdot
x^{3}+1$ where $s=t^{2}$.
\end{tabbing}

\item The Galois groups
$G=G(\Fano)=PSL_3(2)=G(\Serre)=G_{S}$

\item $\forall i,j$ $F_2(P_{i}) \cap
F_2(L_{j})=F_2$

\item The 7 fields $F_2(P_{i})$ correspond to 7
subgroups of
$G$, conjugate and isomorphic to $S_{4}$.

\item The 7 fields $F_2(L_{i})$ correspond to the
other 7 subgroups
of $G$, isomorphic to the previous 7 subgroups
under an outer
automorphism of $G$.

\item For the tower of fields $F_2(t)\subset
F_2(P_{0})\subset
F_2(L_{0},P_{0})\subset R$

\begin{tabular}{ll}
$[K:R]$ & $=3$  (Frobenius automorphism)\\
$[R:F_2(L_{0},P_{0},t)]$ & $=2$\\
$[F_2(L_{0},P_{0},t):F_2(P_{0},t)]$ & $=4$\\
$[F_2(P_{0},t):F_2(t)]$ & $=7$
\end{tabular}

\item $R=F_2(L_{0},P_{0},t,y)$ where $y$ is a zero
of the
polynomial $x^{2}+L_{0}\cdot x+t$.
\end{enumerate}
\end{prop}

Since $\Fano(x)$ is irreducible, and $G$ contains
the Frobenius
automorphism obtained by applying the map
$C\rightarrow C^{2}$ to
the coefficients of the series $P_{i}$, (or
equivalently any series
$S(t)\rightarrow S(t)^{2}$): $7,3$ divide $|G|$
so by
consideration of maximal subgroups $G \cong
PSL_3(2)$ or
$|G|$=$21$.

To show the first isomorphism is correct, show $2\mid \mid G\mid$. Let
\begin{equation*}h(x)=(x-P_{1})\cdot(x-P_{2})\cdot(x-P_{4})=x^{3}+u\cdot
x+L_{0}\end{equation*}
(By the definition of $P_{i}$ and proposition 0.1,
$P_{1}+P_{2}+P_{4}=0$; and $u$,  $L_{0}\in
F_2[[t]]=$ the set of
power series in $t$ with coefficients in $F_2$
(Since the
coefficients of these series are fixed by the
Frobenius
automorphism).

Dividing $h(x)$ into $\Fano(x)$ we get a remainder
\begin{equation*} \label{rx}
r(x)=(t+L_{0}^{2}+u^{3})\cdot x+(1+u^{2}\cdot
L_{0})=0\end{equation*}

Equating the coefficients to zero, we derive the
equation:
\begin{equation}\label{L7}L_{0}^{7}+t^{2}\cdot
L_{0}^{3}+1=0\end{equation}

which shows $L_{0}$ is a zero of $\Serre(x)$.

We also have
\begin{equation}\label{u} u=t\cdot L_{0}+L_{0}^{3}\in
F_2(L_{0},t)\end{equation}

Dividing $\Fano(x)$ by its
factor,$(x-P_{0})\cdot(x^{3}+u \cdot
x+L_{0})$, we get the full factorization

\begin{equation} \label{Fano}
\Fano(x)=(x-P_{0})\cdot(x^{3}+u \cdot
x+L_{0})\cdot(x^{3}+P_{0}\cdot
x^{2}+(u+P_{0}^{2})\cdot x+(L_{0}+u\cdot
P_{0}+P_{0}^{3})\end{equation}

valid over $F_2(P_{0},L_{0},t)$ by (\ref {u})

Again consider: $h(x)=x^{3}=u\cdot
x+L_{0}=(x-P_{1})\cdot(x-P_{2})\cdot(x-P_{4}).$

Since $h(P_{i})=0,i=1,2,4$ for these values of $i$:

\begin{equation*}P_{i}^{3}=u\cdot
P_{i}+L_{0}.\end{equation*}

Let $y=P_{1}^{2}\cdot P_{2}+P_{2}^{2}\cdot
P_{4}+P_{4}^{2}\cdot
P_{1}$ and $y' = P_{2}^{2}\cdot P_{1}+P_{4}^{2}\cdot
P_{2}+P_{1}^{2}\cdot P_{4}$

Then:
\begin{equation} \label{ypy} y+y'=L_{0}\end{equation}
and
\begin{equation} \label {yy} y\cdot y'=t\end{equation}
\begin{quote}
\mbox{\llap{\large[\normalsize\enspace}}\ignorespaces
Let $S(i,j,k)$ be the symmetric polynomial obtained 
by applying the cyclic permutation $(1,2,4)$ to the
subscripts of 
$P_{1}^{i}\cdot P_{2}^{j}\cdot P_{4}^{k}$ 
and adding the products obtained.
\begin{alignat*}{2}
    0 &= S(1,0,0) &&= P_{1}+P_{2}+P_{4} \\
    u &= S(1,1,0) &&= P_{1}\cdot P_{2}+P_{2}\cdot P_{4}+P_{4}\cdot P_{1} \\
L_{0} &= S(1,1,1) &&= P_{1}\cdot
P_{2}\cdot P_{4} 
\end{alignat*}
So: 
\[ 0=S(1,0,0)\cdot S(1,1,1)=y+y'+L_{0}, \]
establishing equation (\ref{ypy}).
\newpars
And: 
\[ y\cdot y'=S(3,3,0)+S(2,2,2)+S(4,1,1). \]
Using (\ref{Fano}):
\begin{align*}
S(3,3,0) &= u^{2}\cdot S(1,1,0)+L_{0}^{2} = u^{3}+L_{0}^{2} \\
S(2,2,2) &= L_{0}^{2} \\
S(4,1,1) &= S(1,1,1)\cdot S(3,0,0) = L_{0}\cdot(u\cdot S(1)+L_{0})
\end{align*}
So: 
\[ y\cdot y'=u^{3}+L_{0}=t \]
by equation (\ref{rx}), establishing equation(\ref{yy}).
\newpars
Since the coefficients of $y$ and $y'$ are fixed by the 
Frobenius automorphism the coefficients are in $F_2$ 
and $y,\,y'\in R$.\mbox{\rlap{\enspace\large]\normalsize}}
\end{quote}
Also, by (\ref{ypy}) and (\ref{yy}),
$L_{0}=\frac{y^{2}+t}{y}$, so
using (\ref{L7}), there is an equation of
degree 16 in
$y$, which decomposes into irreducble factors over
$F_2(t)$:
\[
(y^{14}+t\cdot y^{12}+y^{7}+t^{6}\cdot y^{2}+t^{7})\cdot(y^{2}+t)=0.
\]
So $\mid G \mid$ is even, implying $G \cong
PSL_3(2)$ and
establishing statement (2). By the factorization
(\ref{Fano}) :
$8\nmid[K:F_2(L_{0},P_{0},t)]$

But $[K:F_2(L_{0},t)]$=24. So $P_{0}\notin
F_2(L_{0},t)$

So by (\ref{Fano}) : $(x^{4}+u\cdot x^{2}+L_{0}\cdot
x+u^{2})\cdot(x^{3}+u\cdot x+L_{0})$

are the irreducible factors of $\Fano(x)$ over
$F_2(L_{0},t)$

Thus $[F_2(L_{0},P_{0},t):F_2(P_{0},t)]=4$

and statement(6) follows since the Frobenius
automorphism has order
3.

Now the Fano geometry can be defined by taking points
as elements of
$F_{2^3}^{*}$, and lines as 3-element sets
${\alpha,\beta,\gamma}$
such that $\alpha + \beta + \gamma =0$ . All such sets
are of the
form $\{\alpha=\varepsilon^{i},\beta =
\varepsilon^{j}, \gamma =
\varepsilon^{k}\}$ , where the list
$[i,j,k]\equiv[1,2,4]\pmod7$.
Or, alternatively, $P_{i}+P_{j}+P_{k}=0$ if and only
if
$\{\varepsilon^{i}, \varepsilon^{j},\varepsilon^{k}\}$
is a line,
which establishes statement (1).

The permutation $(1\,2)(3\,6) \in G$ fixes
$L_{0}$ and $P_{0}$
but interchanges $y$ and $y'$ which establishes
statement (7).

Applying $G$ to $F_2(P_{0})\bigcap
F_2(L_{i})=F_2$,
$i=1...6$ establishes statement(3).
Statement (4) follows from Galois
theory and the subgroup structure of $PSL_3(2)$.

For each $i$,$j$ there is an outer automorphism taking the point
stabilizer of $P_{i}$ onto the line stabilizer of $L_{j}$ \cite{Atlas}.
\newpars
This establishes statement(5) and the proof.
\newpars
Finally, there is a natural correspondence
between the series solutions to $s_{11}(x)=x^{11}+tx^{2}-1$ over $F_3(t)$ and
the 11 points of the Steiner system $\Sigma(4,5,11),$ \cite{McWSlo},
which is established by:

\begin{prop}\label{M11}

The polynomial $s_{11}(x) = x^{11}+tx^{2}-1$ has 11 Taylor series solutions
$\theta_{0}(t),\theta_{1}(t),\ldots,\theta_{10}(t)$ where
\begin{enumerate}
\item
The constant term, $\theta_{i}(0)$, is $\varepsilon^i$,
where $\varepsilon$ is
\\
 a primitive $11^{th}$ root of 1 mod 3,
satisfying ${\varepsilon}^{5}+2\,{\varepsilon}^{3}+{\varepsilon}^{2}+2\,\varepsilon+2=0$.

\item

$\theta_{i}=  \sum_{k} c_{k} \varepsilon^{(2k+1)i} t^{k}$
where the $c_{k} \in F_3$ satisfy the recursion relations:
\begin{itemize}
\item[  i.] $c_{3j+2}$ $= 0$
\item[ ii.] $c_{3j+1}$ $= c_{j}$
\item[iii.] $c_{3j}  $ $= \sum_{3m+n=j} c_{m} c_{n} $ mod 3
\end{itemize}

\item

The Galois group of $s_{11}(x)$, $M_{11}$,
is generated by the two permutations,
$\alpha =(0123456789X)$  and $\sigma =(36)(40)(5X)(89)$,
(where $X$ represents the number $10$), 
acting on the subscripts of the $\theta_{i}$.

\item

The $66$ Steiner 5-ads are the orbit of the 5-ad $\{X8267\}$ under the actions of $\alpha$ 
and $\sigma$,
  and to every 5-ad, $\{ijklm\}$, the corresponding product,
$\theta_{i}\theta_{j}\theta_{k}\theta_{l}\theta_{m}$,
  is a zero of the polynomial $f_{66}$ defined by equation(\ref {f66}).

\end{enumerate}

\end{prop}

Let $\theta_{i,0} = \varepsilon^{i}$. The Taylor series solution, $\theta_{i}$, is generated by
the recursion,
\\
$\theta_{i,j+1}= \theta_{i,j}^{12} + t \, \theta_{i,j}^{3}$, establishing (1).
Let this series be $\sum_{k} c_{k} \varepsilon^{\rho(k)} t^{k}$. \\
It follows from the recursion that $c_{k} = 0$ if $k \equiv 2$ mod 3.\\
If $k \equiv 1$ mod 3 and $k = 3j + 1$ then
$c_{k} \varepsilon^{\rho(k)} t^{k} = c_{j}^{3} \varepsilon^{3 \rho(j)i} t^{3j+1} \equiv c_{j}
\varepsilon^{2k+1} t^{k}$ mod 3
if $\rho(j) = 2j+1$,
which establishes (2:i,ii) by induction. \\
If $k \equiv 0$ mod 3, $k=3j$,
since $(a+b)^{12} \equiv \, a^{12} + a^{9}\,b^{3} + a^{3}\,b^{9} + b^{12}$ mod 3,

\begin{align*}
c_{k} \varepsilon^{\rho(k)} t^{k}& = \sum_{9m+3n=k} c_{m}^{9}c_{n}^{3} \varepsilon^{(9 \rho(m) + 3
\rho(n))i} t^{k} \\
 & \equiv  \sum_{9m+3n=k}c_{m} \, c_{n} \, \varepsilon^{(9 \cdot 2m + 3\cdot 2n + 12)i} t^{k}
\textrm{ mod 3}\\
  & = \sum_{3m+n=j}c_{m}\, c_{n} \, \varepsilon^{(2k+1)i} t^{k}
\end{align*}
assuming $\rho(m)= 2m+1$ and $\rho(n)=2n+1$. 
So (2:iii) also follows by induction.
\par\smallskip
For (3) and (4), a calculation shows
that the permutations $\alpha$ and $\sigma$ generate a group of order
$11 \cdot 10 \cdot9 \cdot 8$
and that the orbit of $\{X8267\}$ under this group is the standard set of Steiner 5-ads,
$\Sigma(4,5,11)$,
obtained by the construction in \cite{McWSlo}.
Another computation shows that the corresponding products,
$\theta_{i}\theta_{j}\theta_{k}\theta_{l}\theta_{m}$,
for $\{ijklm\}$  a Steiner 5-ad,
are zeros
of the polynomial $f_{66}$ mod $t^{100}$.
The remaining 396 5-element products have a non-zero term of $t$-degree 4 when substituted
into this polynomial. Since $M_{11}$, the Galois group of $s_{11}(x)$, permutes the zeros of
$f_{66}$,
it is the automorphism group of the set of 5-ads determined by (3), and so is generated by
$\alpha$ and $\sigma$.

\end{document}